\documentclass[11pt]{amsart}

\usepackage{amsmath, amsthm, amssymb}

\usepackage{comment}

\usepackage[abs]{overpic}		  % overpic automatically loads graphicx

\usepackage{float}
\usepackage{wrapfig}
\usepackage{picins}
\usepackage{caption}

\usepackage{xcolor}
\definecolor{indigo}{rgb}{0.29, 0.0, 0.51}  % custom colors
\usepackage[colorlinks, urlcolor=indigo, linkcolor=indigo, citecolor=indigo]{hyperref}

\theoremstyle{plain}

\theoremstyle{definition}

\theoremstyle{remark}

\numberwithin{theorem}{section}

\def\co{\colon\thinspace}

\newcommand{\R}{\mathbb{R}}           % the real numbers
\newcommand{\Z}{\mathbb{Z}}           % the integers
\makeatletter
\newcommand*\bigcdot{\mathpalette\bigcdot@{0.6}}
\newcommand*\bigcdot@[2]{\mathbin{\vcenter{\hbox{\scalebox{#2}{$\m@th#1\bullet$}}}}}
\makeatother

\DeclareMathOperator\tb{tb}                   % Thurston-Bennequin
\DeclareMathOperator\rot{rot}                 % rotation
\DeclareMathOperator{\spins}{\mathfrak{s}}

\DeclareMathOperator{\HFplus}{{{HF}}^+}       % HF+
\DeclareMathOperator{\Dist}{Dist}   
\DeclareMathOperator{\Cont}{Cont}   
\DeclareMathOperator{\Diff}{Diff}   

\begin{document}

% title
\title{Tight and overtwisted contact structures} 

% author information
\author{John B. Etnyre}
\address{Department of Mathematics \\ Georgia Institute of Technology \\ Atlanta \\ Georgia}
\email{etnyre@math.gatech.edu}

%\subjclass[2020]{57R17}

%abstract
%\begin{abstract}
%\end{abstract}

\maketitle
%\tableofcontents

%%%%%%%%%%%%%%%%%%%%%%%%%%%%%%%%%%%%
\section{Introduction}
%%%%%%%%%%%%%%%%%%%%%%%%%%%%%%%%%%%%
The tight versus overtwisted dichotomy has been an essential organizing principle and driving force in 3-dimensional contact geometry since its inception around 1990. In this article, we will discuss the genesis of this dichotomy in Yasha's seminal work and his influential contributions to the theory. We will focus on his classification of overtwisted contact structures, characterization of tightness (and hence the birth of the dichotomy), and work with Thurston allowing for the construction of new tight contact manifolds which have many applications to low-dimensional topology. Along the way, we will discuss further developments by Yasha and others as well as applications of this work outside of contact geometry. We will end with his recent breakthrough with Borman and Murphy defining overtwisted contact structures in higher dimensions. Along the way, we will discuss the further development of the field in the work of Yasha and many others. This discussion is not meant to be comprehensive, but to give a good indication of the depth, breadth, and profound impact of Yasha's work. We will occasionally sketch proofs to give a small indication of what goes into the result (many will be brief sketches based on an amazing class of Yasha's that the author attended in the Fall of 2000 at Stanford), although we refer to the relevant papers for full details of the proofs. 

Throughout this article we will assume the reader is familiar with basic notions from contact geometry, such as the definition of a contact structure and characteristic foliation, what Legendrian and transverse knots are, and what their classical invariants are. All this can be found in many places, such as \cite{Etnyre03, Etnyre05, Geiges08}. (Several of the figures included below were taken from Yasha's original papers.)

\begin{figure}[htb]
\includegraphics[width=0.7\textwidth]{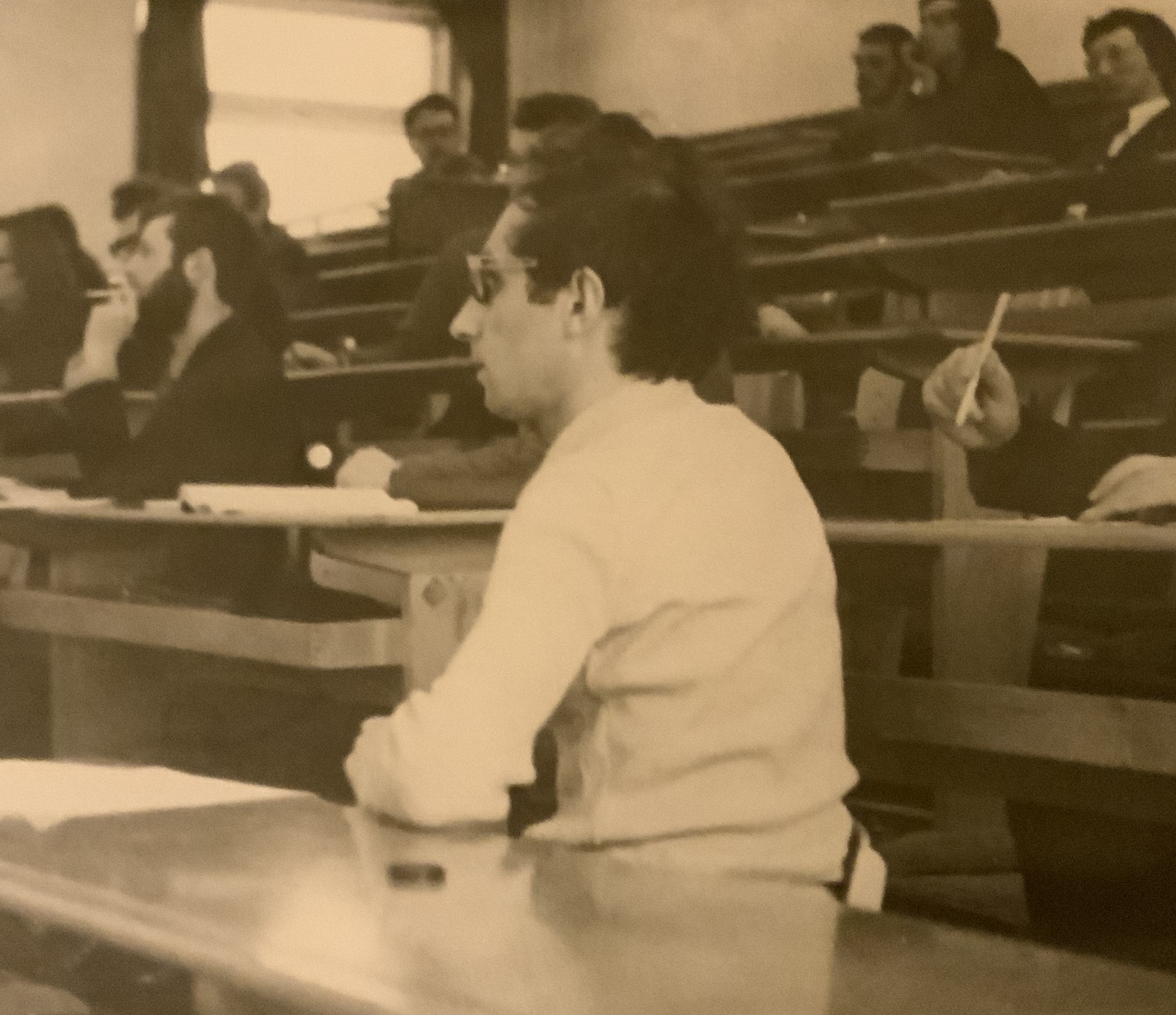}
\caption{Yasha at a conference in the late 1970s when several of the ideas leading to the tight versus overtwisted dichotomy were being developed.}
\end{figure}

%%%%%%%%%%%%%%%%%%%%%%%%%%%%%%%%%%%%
\section{Precursors and the birth of modern contact topology}
%%%%%%%%%%%%%%%%%%%%%%%%%%%%%%%%%%%%

In this section we will confine ourselves to dimension $3$. A fundamental question one can ask is if any oriented $3$-manifold admits a contact structure (note, in dimension $3$ a manifold must be oriented to admit a contact structure). This basic question was answered in 1970 by Martinet \cite{Martinet71} when he showed that indeed any oriented $3$-manifold admits at least one contact structure. The proof relies on the work of Lickorish \cite{Lickorish62} and Wallace \cite{Wallace60} who showed that any oriented $3$-manifold $M$ was obtained from $S^3$ by surgery on some link. To construct a contact structure on the manifold one can realize the link as a transverse link in the standard contact structure $\xi_{std}$ on $S^3$. The contact structure $\xi_{std}$ can be taken to be standard in a neighborhood of each of the components and when performing surgery on the link, the neighborhood is removed and a new solid torus is glued in its place. One can easily extend the contact structure over this glued-in torus to obtain a contact structure on $M$.

The next basic question that arises is whether or not a manifold can support more than one contact structure. A few years after Martinet's work, Lutz \cite{Lutz77} showed that there are contact structures in different homotopy classes of plane fields on a given manifold. He did this by introducing, what is now called, a Lutz twist. 
%As it will be relevant later, we now discuss the Lutz twist. We first consider a standard model of a knot $T$ transverse to a contact structure $\xi$ in a manifold $M$. There will be a standard neighborhood for $T$ that can be described as follows. Consider $\R^2\times S^1$ with the contact structure $\ker(\cos r\, d\phi + r\sin r\, d\theta)$, where $(r,\theta)$ are polar coordinates on $\R^2$ and $\phi$ is the angular coordinate on $S^1$. It is well known that $T$ will have a neighborhood $N$ that is contactomorphic to the neighborhood $N_\epsilon$ in $\R^2\times S^1$, where $N_\epsilon$ is the set of points in $\R^2\times S^1$ with $r$ coordinate less than or equal to $\epsilon$ where we can assume $\epsilon$ is arbitrarily small. We note there is some number $a$ that is in $(\pi, 2\pi)$ such that the characteristic foliation on $\partial N_a$ is the same as on $N_\epsilon$. Given this one can remove $N_\epsilon$ from $M$ and replace it with $N_a$. We note that this does not change $M$, but it does make the contact planes make an extra $\pi$ twist in the solid torus. We say the new contact structure $\xi_T$ is obtained from $\xi$ by a Lutz twist on $T$. 
The key to Lutz's proof is that a Lutz twist can change the homotopy class of the plane field. That is, if $\xi_T$ is obtained from $\xi$ by a Lutz twist then $\xi_T$ and $\xi$ will not be homotopic as plane fields and thus cannot be isotopic as contact structures (and for most choices of $T$ not even contactomorphic). Since it is known that there are infinitely many homotopy classes of plane fields on any given oriented $3$-manifold (see below), we can use the Lutz twist to see that every oriented $3$-manifold admits infinitely many distinct contact structures. 

\pichskip{15pt}% Horizontal gap between picture and text
\parpic[r][b]{%
  \begin{minipage}{0.4\textwidth}
    \includegraphics[width=\linewidth]{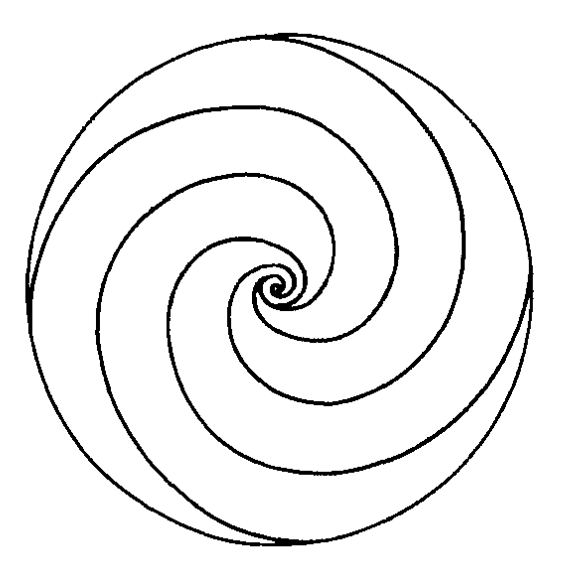}%
    \captionof{figure}{}
    \label{otdisk}
  \end{minipage}
}
We are ready to define an overtwisted contact structure.  
A disk $D$ in a contact manifold $(M,\xi)$ is an overtwisted disk if the characteristic foliation on $D$ has $\partial D$ as a leaf and one singular point inside of $D$. See Figure~\ref{otdisk}. A contact structure $\xi$ on $M$ is called overtwisted if there is an overtwisted disk inside of $M$. The Lutz twist construction above always produces such a disk. 
%, notice that the disk $D=\{(r,\theta,\phi_0): r\leq \pi\}$, for a fixed $\phi_0$, may be slightly perturbed to be an overtwisted disk. 
Thus all the infinitely many distinct contact structures given by Lutz's construction on a given manifold are overtwisted. 

At this point, one might think all contact structures are overtwisted and it might also be true that there is a unique contact structure in each homotopy class of plane field. (We should point out at this point in the historical development of the subject, the definition of overtwisted did not yet exist; so our telling of the story is a bit anachronistic.) While this is not the case, it took the groundbreaking work of Bennequin \cite{Bennequin83} in 1983 to show it was not. 

In \cite{Bennequin83}, Bennequin established the famous Bennequin inequality that says if $T$ is any knot transverse to the standard contact structure $\xi_{std}$ on $\R^3$ (or $S^3$) then its self-linking number satisfies
\[
\mathrm{sl}(T) \leq -\chi(\Sigma)
\]
where $\Sigma$ is any surface with $\partial \Sigma=T$. The self-linking number is the relative Euler class of the contact structure (relative to a natural section of $\xi$ along $T$ coming from the characteristic foliation on $\Sigma$). The only thing we need to know about the self-linking number here is that if one perturbs the boundary of an overtwisted disk appropriately then one will obtain a transverse unknot with self-linking $+1$. Since this clearly violates the Bennequin inequality, the contact structure obtained from a Lutz twist on $S^3$ in the same homotopy class of plane field as $\xi_{std}$ is not isotopic or contactomorphic to $\xi_{std}$. So after Bennequin, we knew that not all contact structures were overtwisted. 

Bennequin's proof of this inequality was ingenious. He took a transverse knot and isotoped it until it was the closure of a braid. He then looked at the singular foliation on the surface $\Sigma$ induced by the braid axis and the constant $\theta$-half planes. (Here we are thinking of $\R^3$ with cylindrical coordinate $(r,\theta, z)$ where the $z$-axis is the braid axis, and a knot is a closed braid if the $\theta$-coordinate restricted to it has no critical points.) Bennequin then isotoped the foliation to simplify it and eventually showed that, for example, if $\Sigma$ was a disk then one could arrange, by possibly changing $T$, but only so that $\mathrm{sl}(T)$ increased, that the foliation on $\Sigma$ contained one singular point and radial lines on the disk. From here one can easily conclude that $\mathrm{sl}(T)=-1$ and of course $\chi(\Sigma)=1$. 

%One indication the central nature of the Bennequin inequality is that there are many ways to prove it ranging from topological, to geometric, and analytic. 

This result is really the birth of modern contact topology and the first hint of many subtle connections between contact geometry and the topology of $3$-manifolds. More specifically, it is at this point that we finally know that contact structures are not just determined by their homotopy data so there is something deeper about their structure than ``just algebraic topology". In addition, we see that the Bennequin inequality shows that one may use contact geometry to bound $-\chi(\Sigma)$ for any Seifert surface for $T$ and since $-\chi(\Sigma)=2g(\Sigma)-1$ (where $g(\Sigma)$ is the genus of $\Sigma)$ we see that contact geometry can give a lower bound on the genus of a Seifert surface for a knot. Determining the minimal genus is a difficult problem that only recently has a fairly tractable solution. 

Overtwisted contact structures (anachronistically) appeared in a few other works, such as \cite{Erlandsson} and \cite{GonzaloVarela83}, but apart from Bennequin's work, where the existence of an overtwisted disk was key, the overtwisted disk seemed more a byproduct of a construction rather than a key feature. 

We now come to the seminal work of Yasha where he first defines what it means for a contact structure to be overtwisted and then completely classifies overtwisted contact structures. This occurred in his 1989 paper \cite{Eliashberg89}. We will discuss this result in detail below, but now continue the journey to the tight versus overtwisted dichotomy. 

To this end, Yasha noticed that a generalization of the Bennequin inequality for close surfaces, which he established in \cite{Eliashberg92a}, was very reminiscent of inequalities for taut foliations. Yasha thought the word ``taut" was overused so decided to call such contact structures tight. Specifically, in 1992 Yasha  \cite{Eliashberg92a} defined a tight contact structure to be one that contained no disks whose characteristic foliation contained a limit cycle. In that same paper, he proved that if a contact structure was not tight it was overtwisted. {\em Thus the birth of the tight versus overtwisted dichotomy.} In his paper  \cite{Eliashberg93}, Yasha also shows a contact structure is tight if and only if the Bennequin inequality is true and, in the paper \cite{Eliashberg92a}, that $S^3$ has a unique tight contact structure (thus Yasha completely classified all contact structures on $S^3$).  

It is important to appreciate Yasha's classification of contact structures on $S^3$. After Bennequin's work, it was thought that this was just the first example of a ``zoo" of non-equivalent contact structures on $S^3$ in the given homotopy class of plane field. Quoting from Yasha's paper \cite{Eliashberg92a}, ``Twenty years ago Jean Martinet (see \cite{Martinet71}) showed that any orientable
closed $3$-manifold admits a contact structure. Three years later after the
work of R. Lutz (see \cite{Lutz77}) and in the wake of the triumph of Gromov's $h$-principle, it seemed that the classification of closed contact 3-manifolds was
at hand. Ten years later in the seminal work \cite{Bennequin83}, D. Bennequin showed
that the situation is much more complicated and that the classification of
contact structures on $3$-manifolds, and even on $S^3$, was not likely to be
achieved." Yasha's work in \cite{Eliashberg89, Eliashberg92a} showed that (1) any homotopy class of plane field has a unique overtwisted contact structure and (2) there is exactly one homotopy class that has another contact structure. So the ``zoo of contact structures" on $S^3$ turns out to be quite well-behaved!

We briefly sketch Yasha's proof that $S^3$ has a unique tight contact structure (up to isotopy). Since Darboux's theorem says that contact structures are all the same in a neighborhood of a point, this result follows by showing that two tight contact structures on the $3$-ball $B^3$ that induce the same characteristic foliation on the boundary are actually isotopic. That is, we need to see that given any tight contact structure $\xi$ on $B^3$ then there is a fixed tight contact structure $\xi_{m}$ that is isotopic to $\xi$. To do this, Yasha gives an ingenious construction of ``model contact structures on the ball" given a characteristic foliation on its boundary that could come from a tight contact structure. That is, given a singular foliation $\mathcal{F}$ on $\partial B^3$ he constructs a model contact structure $\eta_F$ on $B^3$ that induces $\mathcal{F}$ on boundary. Now given $\xi$ on the unit ball $B^3$ we can let $\mathcal{F}_t$ be the characteristic foliation on the boundary of the ball $B_t$ of radius $t$ in $B^3$ and let $\eta_t$ be the model contact structure on the ball whose boundary has characteristic foliation $\mathcal{F}_t$. Now we can define contact structure $\xi_t$ to be $\eta_t$ on $B_t$ and $\xi$ on $B^3-B_t$. Notice that this is a one-parameter family of contact structures on $B^3$ that starts with $\xi$ and ends with $\eta_1$. That is $\xi$ is isotopic to the model contact structure $\eta_1$. (Note there is an issue with $\eta_t$ when $t$ approaches $0$, but we can use Darboux's theorem again to make sure $\xi$ is standard on $B_\epsilon$ and then run the argument from there). Clearly the key to this proof is the construction of the model contact structures, but for that we refer the reader to Yasha's paper \cite{Eliashberg92a}. 

A natural question now is: are there any more tight contact structures? In \cite{Eliashberg90b} Yasha proved that any symplectically fillable contact structure is tight by showing that transverse knots in this contact structure satisfy the Bennequin inequality (which by an earlier result of Yasha implies that it must be tight). The proof involves the analysis of pseudo-holomorphic disks in the symplectic filling. See Geiges' article {\em Filling by holomorphic disks} in this Celebratio. This result gives a large family of tight contact manifolds and there are lots of constructions of symplectically fillable contact structures, see for example \cite{Eliashberg90a, Gompf98}. Later we will discuss a technique developed by Yasha and Thurston to build tight contact structures by ``perturbing taut foliations". Though these contact structures are still symplectically fillable, previous constructions of symplectic fillings could not approach what can be done with this construction.

%%%%%%%%%%%%%%%%%%%%%%%%%%%%%%%%%%%%
\section{The classification of overtwisted contact $3$-manifolds}
%%%%%%%%%%%%%%%%%%%%%%%%%%%%%%%%%%%%
We begin by stating the simplest version of Yasha's classification of overtwisted contact structures on $3$-manifolds. To this end, we let $\Dist(M)$ be the space plane fields on $M$ and $\Cont_{ot}(M)$ be the space of overtwisted contact structures on $M$. Clearly there is an inclusion of $\Cont_{ot}(M)$ into $\Dist(M)$.  Yasha's main result in \cite{Eliashberg89} is that this inclusion induces a bijection
\[
i_*\co \pi_0(\Cont_{ot}(M))\to \pi_0 (\Dist(M)).
\]
That is every homotopy class of plane field contains a unique isotopy class of overtwisted contact structure (recall that contact structures that are homotopic through contact structures are isotopic). 

Understanding homotopy classes of plane fields is fairly straightforward and so Yasha's work gives an effective and complete classification of overtwisted contact structures. We now elaborate on homotopy classes of plane fields. One can show, see \cite{Gompf98}, that there is a map (that requires choices) 
\[
 f\co \pi_0 (\Dist(M))\to H_1(M)
\]
and that $f^{-1}(h)=\Z/(2d(h))$ where $d(h)$ is the divisibility of $h$ in $H_1(M)$. Moreover, $f(\xi)$ is determined by the spin$^c$ structure induced by $\xi$ and there are explicit formulas for distinguishing elements of $f^{-1}(h)=\Z/(2d(h))$, see \cite{Gompf98}. So we see, among other things, that every spin$^c$ structure admits several overtwisted contact structures but also that any oriented $3$-manifold admits infinitely many such structures (because $d(0)=0$). Just for some context, we note that for a lens space $L(p,q)$ the homotopy classes of plane fields are in one-to-one correspondence with 
\[
H_1(L(p,q)\oplus \Z\cong \Z/p\Z\oplus \Z
\]
and on $S^1\times S^2$ the homotopy classes of plane fields are in one-to-one correspondence with 
\[
\oplus_{n\in \Z} (\Z/(2n)\Z)\cong \Z\oplus (\Z/2\Z\oplus \Z/2\Z)\oplus (\Z/4\Z\oplus \Z/4\Z)\oplus \cdots.
\]

We note that this is a classical example of an $h$-principle. Recall, very roughly, a weak form of an $h$-principle is that if there is no topological obstruction for something being true then it is true, and if two things look the same on the topological level then they are the same. Yasha is responsible for many other such $h$-principles, but this is certainly a very influential one, and he has an even stronger version. 

To state the stronger version we denote by $\Dist_p(M)$ the space of plane fields on $M$ that agree with a fixed plane at $p$ and $\Cont^D_{ot}(M)$ the space of contact structures on $M$ that have a fixed disk $D\subset M$ as a standard overtwisted disk that is tangent to $\xi$ at $p$ (by standard overtwisted disk we just mean to fix one characteristic foliation on $D$ as described in the definition of overtwisted disk). Yasha's main theorem in \cite{Eliashberg89} is that the inclusion
\[
\Cont^D_{ot}(M) \to \Dist_p(M)
\]
is a weak homotopy equivalence. It is a simple exercise to check that this result implies the result above about $\pi_0$ follows (note this is not automatic since all the spaces involved are different). 

 \pichskip{15pt}% Horizontal gap between picture and text
\parpic[r][b]{%
  \begin{minipage}{0.4\textwidth}
    \includegraphics[width=\linewidth]{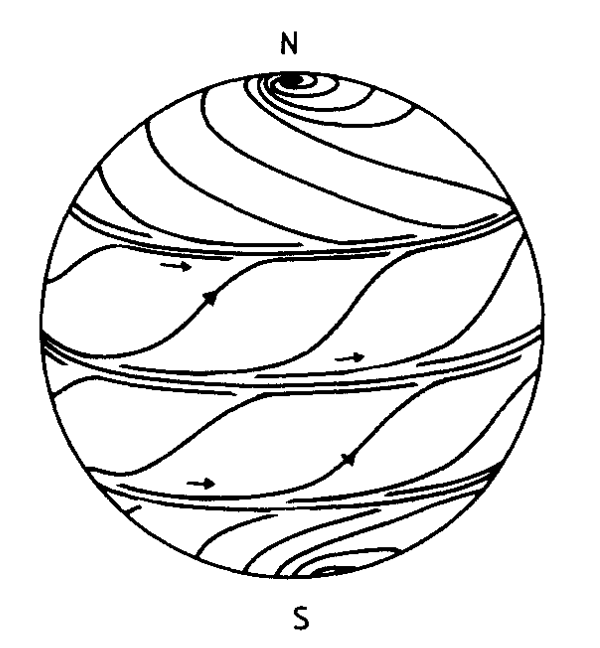}%
    \captionof{figure}{}
    \label{ahorizontal}
  \end{minipage}}
We briefly discuss Yasha's proof that there is a bijection between the homotopy classes of plane fields on a $3$-manifold $M$ and overtwisted contact structures on the manifold. The surjectivity of $i_*$ above follows from the work of Lutz mentioned above, so we consider two overtwisted contact structures $\xi_0$ and $\xi_1$ on a manifold $M$ that are homotopic as plane fields. It is easy to isotope one of them so that they both share an overtwisted disk, say $D$, and hence they agree near $D$. Let $\xi_t, t\in [0,1],$ be the homotopy between the two contact structures. 
One now picks a triangulation of $M$, then one may homotop all the $\xi_t$, fixing $\xi_0$ and $\xi_1$, so that they agree in a neighborhood of $D$ and near the vertices of the triangulation.  One may then further homotop the $\xi_t$, fixing $\xi_0$ and $\xi_1$, so that they are all contact structures on the $2$-skeleton of the triangulation. Now the $\xi_t$ are contact structures except possibly on a finite union of balls $B_i$.
If the triangulation is chosen sufficiently small (with respect to some auxiliary metric) then one can arrange that the characteristic foliation on the boundaries of the $B_i$ are almost horizontal, see Figure~\ref{ahorizontal}.
 One can now find transverse curves connecting the singularities of one of the $B_i$ to another (and also to the neighborhood of $D$). A neighborhood of the balls and the curves is a single ball and the $\xi_t$ have been homotoped to be contact outside this one ball. Now one can use the fact that we have an overtwisted disk to build a model contact structure to fill in these balls. This provides the isotopy of contact structures from $\xi_0$ to $\xi_1$.

%%%%%%%%%%%%%%%%%%%%%%%%%%%%%%%%%%%%
\section{Legacy and applications of overtwisted contact structures}
%%%%%%%%%%%%%%%%%%%%%%%%%%%%%%%%%%%%
It was common after Yasha's classification of overtwisted contact structures for people (including myself) to say that overtwisted contact structures were not interesting since they were determined by their homotopy class of plane field whereas tight contact structures were seen to be tightly connected to the topology of $3$- and $4$-manifolds. But just because something is classified, and fairly easy to understand, does not make it uninteresting! For example, surfaces have been classified for a long time now, but the study of surfaces is still a robust field of mathematics. One just needs to ask more refined questions, such as what do we know about diffeomorphisms of surfaces or how can we apply them to understand $3$-manifolds, among many other questions. The last few decades have shown that overtwisted contact structures are similarly interesting. There are many interesting and subtle features of overtwisted contact structures as well as several exciting applications. We will discuss these now.

\subsection{Results about overtwisted contact structures}
We begin with Legendrian knots in an overtwisted contact structure. A Legendrian knot $L$ in an overtwisted contact structure is called loose if the contact structure restricted to its complement is also overtwisted; otherwise, the knot is called non-loose (the term exceptional is also used). It is a folklore result, see \cite{Etnyre13}, that any two loose knots with the same classical invariants (that is, they are in the same knot type, have the same Thurston-Bennequin invariant, and have the same rotation number) are coarsely equivalent, which means that there is a contactomorphism of the ambient manifold that is smoothly isotopic to the identity and takes one of the knots to the other. We will see, below, that understanding loose knots up to isotopy through Legendrian knots is much more subtle. But first, we turn to non-loose knots. 

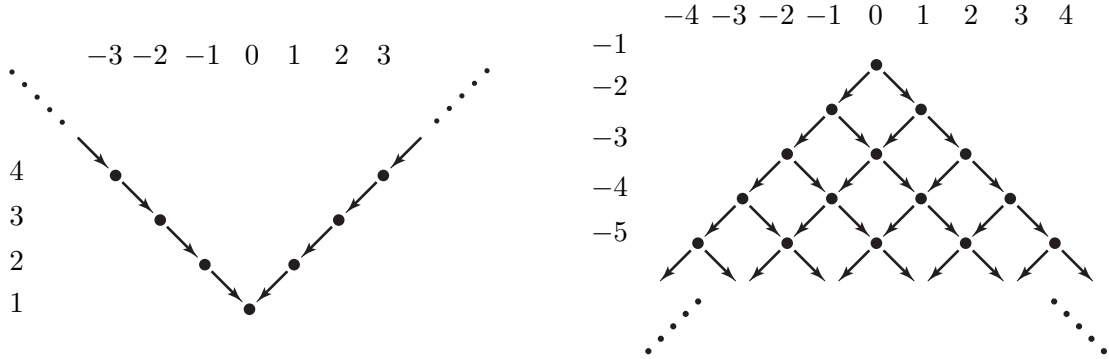
\begin{figure}[htb]{
\begin{overpic}%[grid,tics=10] 
{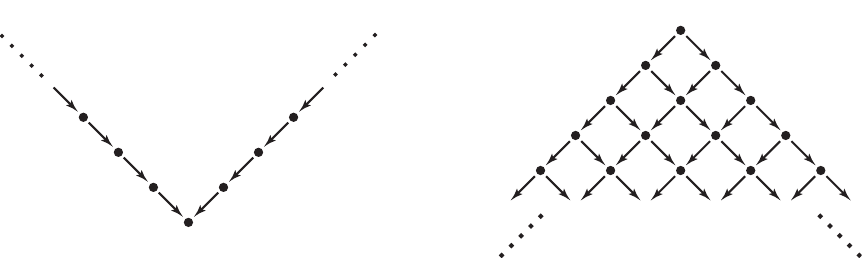}
\put(89, 110){$0$}
\put(66, 110){$-1$}
\put(46, 110){$-2$}
\put(29, 110){$-3$}
\put(105, 110){$1$}
\put(123, 110){$2$}
\put(139, 110){$3$}
\put(0, 16){$1$}
\put(0, 32){$2$}
\put(0, 49){$3$}
\put(0, 66){$4$}
%%%
\put(325, 125){$0$}%90+235
\put(343, 125){$1$}
\put(361, 125){$2$}
\put(380, 125){$3$}
\put(397, 125){$4$}
\put(301, 125){$-1$}
\put(283, 125){$-2$}
\put(265, 125){$-3$}
\put(247, 125){$-4$}
\put(220, 114){$-1$}
\put(220, 98){$-2$}
\put(220, 79){$-3$}
\put(220, 60){$-4$}
\put(220, 43){$-5$}
\end{overpic}}
\caption{On the left, each dot represents a non-loose Legendrian unknot and the arrows indicate what happens after stabilization. On the left we see the same for Legendrian unknots in the standard tight contact structure on $S^3$.}
\label{fig:unknotmr2}
\end{figure}

In \cite{EliashbergFraser09} Yasha and Fraser classified non-loose Legendrian unknots in $S^3$. We recall that overtwisted contact structures on $S^3$ can be indexed by the integers. So we let $\xi_n$ be the overtwisted contact structure on $S^3$ with Hopf invariant $n$. Yasha and Fraser's classification says that non-loose unknots exist only in $\xi_{-1}$ and in $\xi_{-1}$ the complete list, up to coarse equivalence, is $L_1$ and $L^\pm_n$ for $n\geq 2$, where 
\[
\tb(L_1)=1, \tb(L_n^\pm)=n, \rot(L_1)=0, \text{ and } \rot(L_n^\pm)=\pm n-1,
\]
and the $\mp$-stablization of $L_n^\pm$ is $L_{n-1}^\pm$ for $n>2$ and the $\mp$-stablization of $L_2^\pm$ is $L_1$. It is easiest to visualize this classification by plotting the rotation number and Thurston-Bennquin invariants of these non-loose knots; such a plot is called the ``mountain range" of the non-loose Legendrian unknots. See Figure~\ref{fig:unknotmr2}. 
It is interesting to contrast this with the classification of Legendrian unknots in the tight contact structure on $S^3$. This was given in \cite{EliashbergFraser09} too and the mountain range for these knots is shown in Figure~\ref{fig:unknotmr2}. We note that in a tight contact manifold, there is a Legendrian version of the Bennquin inequality which says for a Legendrian knot $L$ that bounds a surface $\Sigma$ we have
\[
\tb(L)+\rot(L)\leq -\chi(\Sigma). 
\]
This inequality shows that in a tight contact manifold, the mountain range for any knot is bounded above, but we see in the classification of non-loose unknots that they can have Thurston-Bennequin invariant arbitrarily large! This seems quite surprising, but it also seems to be a general phenomenon. For example, in work of Min, Mukherjee, and the author \cite{EtnyreMinMukherjee22Pre} they classify all non-loose Legendrian torus knots in overtwisted contact structures on $S^3$ up to coarse equivalence, and the mountain range for any torus knot is not bounded above! They also show that for any fixed torus knot there are only a finite number of overtwisted contact structures that admit non-loose representatives of that knot. So we see when studying non-loose knots there are many subtle features and there are many open questions about them. For example, given a knot $K$ in $S^3$ is it always true that $K$ admits non-loose representatives in only finitely many overtwisted contact structures on $S^3$? And how does one determine which contact structures those are? Are the mountain ranges for non-loose knots always unbounded from above? Why?

Now turning to the classification of non-loose unknots up to Legendrian isotopy we see the beautiful work of Vogel \cite{Vogel2018}. He showed that given a pair of integers $(r,t)$ with $r+t$ odd, there is a unique, up to Legendrian isotopy, loose Legendrian knot $L$ in $(S^3,\xi_{-1})$ with $\tb(L)=t$ and $\rot(L)=r$ if $t<0$, but if $t>0$ then there are two distinct loose Legendrian unknots in $(S^3,\xi_{-1})$ with these invariants (though of course from above, they will be coarsely equivalent). This is quite surprising! Moreover, he showed that in $(S^3,\xi_n)$ with $n\not=-1$, loose Legendrian unknots are classified, up to Legendrian isotopy, by their rotation number and Thurston-Bennequin invariant. Turning to non-loose Legendrian unknots, Vogel showed that for each of the Legendrian knots in Yasha and Fraser's classification, there are two Legendrian representatives up to isotopy. Again, very surprising and points to a subtle difference between the classification of Legendrian knots in overtwisted contact structure up to Legendrian isotopy and coarse equivalence.

We now turn to the contactomorphisms of overtwisted contact structures. We will denote by $\Diff_+(M,\xi)$ the space of contactomorphisms of $\xi$ that preserve the orientation on $\xi$.  In \cite{Vogel2018} Vogel showed that
\[
\pi_0(\Diff_+(S^3,\xi_n))=\begin{cases}
\Z/2\Z\oplus \Z/2\Z & \text{ if } n=-1\\
\Z/2\Z& \text{ otherwise.}
\end{cases}
\]
This result was originally announced by Chekanov, see \cite{EliashbergFraser09}, but was first proven in \cite{Vogel2018}. 
We note that the $\Z/2\Z$ factor in all of the diffeomorphism groups is detected by invariant due to Dymara \cite{Dymara01} while the other $\Z/2\Z$ factor for $\xi_{-1}$  is detected by the action of a contactomorphism on the non-loose unknots. 
We contrast the result above with Yasha's proof \cite{Eliashberg92a} that
\[
\pi_0(\Diff_+(S^3,\xi_{tight}))=\{1\},
\]
that shows contactomorphisms of overtwisted contact structures on $S^3$ is ``more interesting" than the contactomorphisms of the tight contact structure!

\subsection{Applications of overtwisted contact structures}
We will focus on three applications of overtwisted contact structures: the surgery description of all contact $3$-manifolds, the existence of near symplectic forms on $4$-manifolds with positive $b_2^+$, and the existence of Engel structures on parallelizable $4$-manifolds. 

 \pichskip{15pt}% Horizontal gap between picture and text
\parpic[r][b]{%
  \begin{minipage}{0.4\textwidth}
    \includegraphics[width=\linewidth]{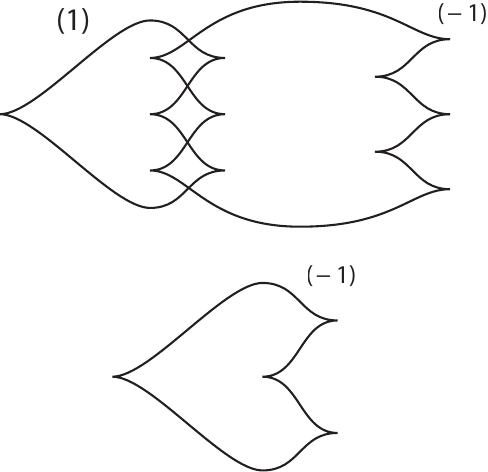}%
    \captionof{figure}{}
    \label{2ex}
  \end{minipage}}
One of the first applications of overtwisted contact structures is to contact geometry in general. In \cite{DingGeiges04} Ding and Geiges showed that any contact $3$-manifold $(M,\xi)$ could be obtained from $S^3$ with its standard tight contact structure by contact surgery on a Legendrian link. For example, Figure~\ref{2ex} shows surgery diagrams for two different overtwisted contact structures on $S^3$. It is a beautiful result that does not mention overtwisted contact structures but essentially uses them in its proof.  Moreover, it should be noted that this result is now ubiquitous in contact geometry; it gives us one of the key ways to present general contact manifolds! The proof of this result goes as follows. It is easy to show that there are Legendrian knots in any contact manifold such that $+1$ contact surgery on the knot yields an overtwisted manifold. So given $(M,\xi)$ we can find a Legendrian knot $L$ on which we can do $+1$ contact surgery to get a manifold $(M',\xi')$ with an overtwisted contact structure. We note that in $(M',\xi')$ there is ``dual" Legendrian $L'$ on which we can perform $-1$ contact surgery to obtain $(M,\xi)$ back again. Now there is also a Legendrian $K$ in $(S^3,\xi_{tight})$ on which we can perform $+1$ contact surgery to obtain $(M'',\xi'')$ where $\xi''$ is overtwisted. It is now fairly easy to see that one can get from any overtwisted contact manifold to any other overtwisted contact manifold by a sequence of $-1$ contact surgeries, see for example \cite{EtnyreHonda02a}. This completes the proof.

% \pichskip{15pt}% Horizontal gap between picture and text
%\parpic[l][b]{%
%  \begin{minipage}{0.4\textwidth}
%    \includegraphics[width=\linewidth]{fig/sympsing}%
%    \captionof{figure}{}
%    \label{simpsing}
%  \end{minipage}}
There are also several applications to the geometry of $4$-manifolds. We first mention a result of Gay and Kirby \cite{GayKirby04} that says a $4$-manifold having $b_2^+>0$ admits a near symplectic structure. A near symplectic structure on a $4$-manifold is a closed $4$-form that is a symplectic form away from some embedded circles and along those circles it has a special form. The key to Gay and Kirby's construction is to realize that near the singular circles, one could build a model for the near symplectic structure whose boundary is overtwisted. They then use some clever handle manipulations to extend this structure over the rest of the manifold and to do this, it is essential to know that overtwisted contact structures are determined by their homotopy class of plane field. The reader can refer to Gay's article in "Celebratio Mathematica for Kirby" for a discussion of how this near symplectic construction led to their discovery of trisection of $4$-manifolds which is currently a very active area of research.

Another application of overtwisted contact structures is the existence of Engel structures on parallelizable $4$-manifolds. An Engel structure on a $4$-manifold $X$ is a $2$-plane field $\mathcal{D}$ on $X$ such that $[\mathcal{D},\mathcal{D}]$ has rank $3$ everywhere and $[\mathcal{D},[\mathcal{D},\mathcal{D}]]=TX$. Here $[\mathcal{D},\mathcal{D}]$ is the subbundle of $TM$ formed by taking all the commutators of local sections of $\mathcal{D}$. Engel structures are interesting for many reasons, but one particularly interesting feature is that they are stable in the sense that a $C^2$ small perturbation of an Engel structure is an Engel structure. What makes this remarkable is that Mongomery \cite{Montgomery93} classified all stable distributions, and they come in three broad families (non-singular line fields, contact structures, and ``even" contact structures) and Engel structures on $4$-manifolds! It is not hard to check that if a $4$-manifold admits an Engel structure then it must be parallelizable. There has been some work trying to prove that all parallelizable $4$-manifolds admit Engel structures. In 2009, Vogel \cite{Vogel09} finally established this. The key idea of Vogel was to take a ``round handle" decomposition of the $4$-manifold, build model Engel structures on the handles and try to glue the Engel structures together. To achieve the gluing, it is essential that the boundaries of the round handles have overtwisted contact structures induced on them. It is the flexibility that the overtwistedness allows for which makes the construction work.

%%%%%%%%%%%%%%%%%%%%%%%%%%%%%%%%%%%%
\section{Legacy and applications of tight contact structures}
%%%%%%%%%%%%%%%%%%%%%%%%%%%%%%%%%%%%
From the very beginning, tight contact structures seemed more interesting than overtwisted ones since their existence is more subtle, and we know that the Bennequin inequality is true for tight contact structures and thus there is a strong connection between such contact geometry and the topology of the manifold supporting the structure. In this section, we will discuss how our knowledge of tight structures has progressed since its inception, including an important result of Yasha and Thurston, the current state of our understanding of tight contact structures, and applications of tight contact structures.

%%%%%%%%%%%%%%%%%%%%%%%%%%%%%%%%%%%%
\subsection{The Eliashberg-Thurston theorem}
%%%%%%%%%%%%%%%%%%%%%%%%%%%%%%%%%%%%
We mentioned above that Yasha showed symplectically fillable contact structures had to be tight. This allows us to construct many tight contact manifolds, but in 1996, Yasha and Thurston introduced a powerful new construction of tight contact manifolds that has been a key part of many major advances in $3$-manifold topology. 

A foliation on a $3$-manifold $M$ is a ``nice" decomposition of $M$ into a union of surfaces. If you look at the tangent planes to these surfaces you get a distribution $\mathcal{F}$ that determines the foliation. Yasha and Thurston's remarkable theorem \cite{EliashbergThurston96, EliashbergThurston98} said that any oriented $C^2$ foliation of an oriented $3$-manifold $M$, other than the foliation of $S^1\times S^2$ by $\{pt\}\times S^2$, can be $C^0$ approximated by a positive $\xi_+$ and negative $\xi_-$ contact structure. 

While this is a very pleasing result in its own right, further work of Yasha and Thurston made it a key tool in low-dimensional topology. Specifically, if the foliation $\mathcal{F}$ on $M$ were taut, then they could build a symplectic structure on $M\times [-1,1]$ that gave a filling of $(M\times\{1\}, \xi_+)\cup (-M\times\{-1\},\xi_-)$. By taut, we mean that on $M$ there is a volume-preserving vector field transverse to $\mathcal{F}$. Of course this means that $(M,\xi_+)$ and $(-M,\xi_-)$ are symplectically fillable and hence tight. 

This last result is particularly powerful given work of Gabai, \cite{Gabai83}. Gabai showed that any closed, oriented, irreducible $3$-manifold $M$ with $H_2(M)$ non-trivial admits a taut foliation. Combined with Yasha and Thurston's work, we know that such a manifold always admits a tight contact structure! We now know that, in some sense, most oriented, irreducible $3$-manifolds admit tight contact structures; in fact, the only way an irreducible $3$-manifold might not admit a tight contact structure is if it is a rational homology sphere. 

To prove their result, Yasha and Thurston introduced the notion of a confolation. This is a plane field $\xi$ defined by a $1$-from $\alpha$ such that $\alpha\wedge d\alpha\geq 0$. So there are regions where $\xi$ defines a foliation, and regions where it defines a contact structure. The proof that a foliation (other than the standard foliation of $S^1\times S^2$) can be perturbed to a contact structure has two steps. In Step 1, one shows that given a foliation $\xi$ on $M$ one may $C^0$-perturb it to a confoliation $\xi'$ such that every point in $M$ can be connected by a path tangent to $\xi'$ to a region where $\xi'$ is a contact structure. In Step 2, one shows that given a confolation as in Step 1, one can $C^{\infty}$ deform it into a contact structure $\xi$. 

 \pichskip{15pt}% Horizontal gap between picture and text
\parpic[r][b]{%
  \begin{minipage}{0.5\textwidth}
    \includegraphics[width=\linewidth]{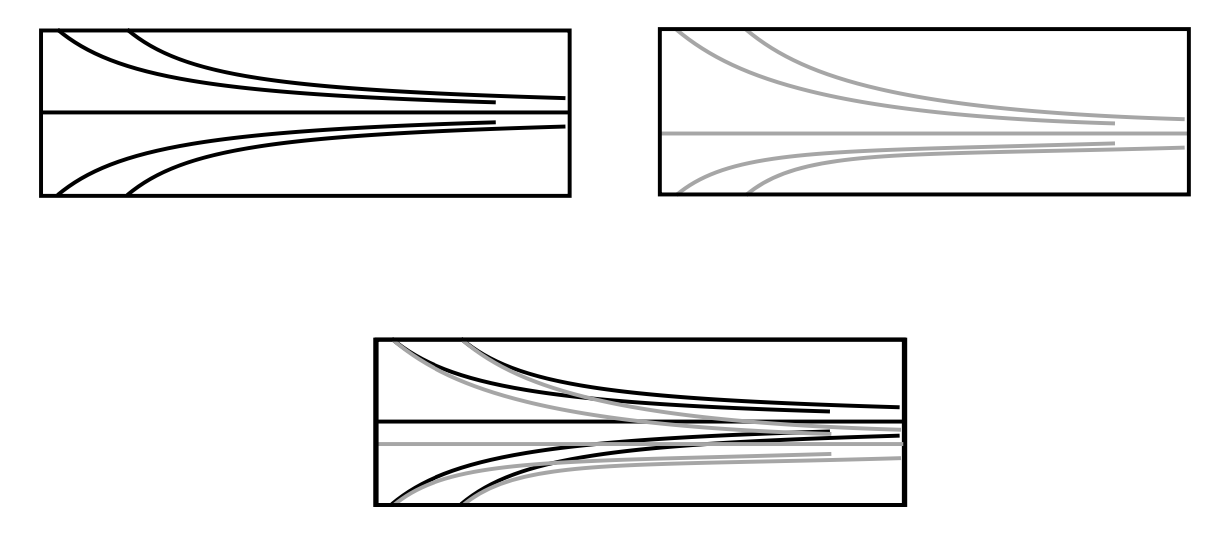}%
    \captionof{figure}{}
    \label{sheer}
  \end{minipage}}
  The key idea for Step 1 is that if one has a leaf $L$ of the foliation $\xi$ that contains a closed curve $\gamma$ with non-trivial linear holonomy, then one may ``shear" $\xi$ in a neighborhood of $\gamma$ to get a new plane field that is contact in a neighborhood of $\gamma$ (and unchanged outside of a slightly larger neighborhood). To describe this more fully, we recall the definition of holonomy. If $A$ is an annulus containing $\gamma$ that is transverse to $\xi$ then there is a non-singular $1$-dimensional foliation of $A$ induced from $\xi$, and $\gamma$ is a leaf in this foliation. Taking a transverse curve $\eta$ to the foliation on $A$, we can define a first return map of the line field on $A$ by starting with a point on $\eta$ and pushing it around the leaf of the foliation on $A$ until we return to $\eta$. We can assume that $\eta=[-1,1]$ and $0$ is the intersection of $\gamma$ with $\eta$. Thus $0$ is a fixed point of the first return map. We call the first return map the holonomy of $\gamma$. Assuming the holonomy has derivative at $0$ that is not $1$, then we say it has nontrivial linear holonomy. Consider a neighborhood $[-\epsilon,\epsilon]\times A$ of $A$. If we shift the foliation up slightly along $A\times \{\epsilon\}$, then we can replace $\xi$ on the interior of $[-\epsilon,\epsilon]\times A$ with a contact structure. In Figure~\ref{sheer}, we see the $A\times\{-\epsilon\}$ and $A\times\{\epsilon\}$ (after sheering) on the top, and on the bottom they are superimposed. From this, one can see a plane field that is tangent to the $[-\epsilon,\epsilon]$ direction and rotates from the foliation on one annulus to the other will be a contact structure. Yasha and Thurston also show that given other non-trivial holonomies for $\gamma$ one may still $C^0$ perturb $\xi$ in a similar manner (though it is more complicated). 

To complete Step 1, Yasha and Thurston show that one can find sufficiently many curves with non-trivial holonomy. To this end, we recall that a subset of $M$ is called minimal if it is a non-empty closed set that is a union of leaves and does not contain a smaller such set. It is easy to see that every leaf in a foliation limits to a minimal set. So if we can find the appropriate curves $\gamma$ with good holonomy in the minimal sets, then we can easily prove that perturbing as above will produce the desired $\xi'$. Luckily, we know that a minimal set in $M$ must be either (1) all of $M$, (2) a closed leaf, or (3) an ``exceptional" minimal set. One can then analyze these cases to find the desired $\gamma$ (possibly after altering the foliation!). 

There are two approaches to Step 2. For the first approach, Yasha and Thurston show how to take curves $\nu$ tangent to $\xi'$ that starts in a contact region of $M$ and ends at a foliated region of $M$, and then perturb $\xi'$ in a $C^\infty$ way to make it contact in a neighborhood of $\nu$. One may iterate this process to turn $\xi'$ into a contact structure. In the second approach, one may apply a result of Altshuler \cite{Altschuler95} to deform $\xi'$ into a contact structure. Altschuler defines a kind of ``heat" flow on the space of differential $1$-forms that can be used to ``distribute" the twisting of the plane field across the manifold.

%%%%%%%%%%%%%%%%%%%%%%%%%%%%%%%%%%%%
\subsection{The current state of tight contact structures}
%%%%%%%%%%%%%%%%%%%%%%%%%%%%%%%%%%%%
The most basic question about tight contact structures is when they exist. As mentioned in the last sub-section, the work of Yasha, Thurston, and Gabai shows that the only way an irreducible $3$-manifold might not admit a tight contact structure is if it is an integral homology sphere. It turns out that there are manifolds that do not admit a tight contact structure.  The first such manifold was found in 2001 in the work of Honda and the author \cite{EtnyreHonda01a} who show that the Poincar\'e homology sphere with its ``non-standard" orientation does not admit a tight contact structure, but with its standard orientation, it does admit a tight contact structure. It is currently unknown if every irreducible integer homology sphere admits a tight contact structure with at least one orientation. 

Generalizing the Poincar\'e homology sphere example, Lisca and Stipsicz \cite{LiscaStipsicz07} showed that $2n-1$ surgery on the $(2,2n+1)$-torus knot does not admit any tight contact structures (notice that the Poincar\'e homology sphere with its non-standard orientation is the $n=1$ case). They later showed that a Seifert fibered space admits a tight contact structure if and only if it is not $2n-1$ surgery on the $(2,2n+1)$-torus knot \cite{LiscaStipsicz09}. Thus, understanding the existence of tight contact structures on a given irreducible $3$-manifold reduces understanding tight contact structures on hyperbolic rational homology spheres. Here virtually nothing is known.%, though see \cite{ConwayMin20} 

Knowing something about the existence of tight contact structures on a given manifold, one might ask how many tight contact structures can exist on a given $3$-manifold. The first result along these lines was again due to Yasha. In \cite{Eliashberg92a}, Yasha showed that if $\xi$ was a tight contact structure and $\Sigma$ was any embedded surface, then we have an adjunction inequality 
\[
|\langle e(\xi),[\Sigma]\rangle| \leq 2g(\Sigma)-2
\]
where $\langle e(\xi),[\Sigma]\rangle$ is the pairing of the Euler class of $\xi$ with the homology class of $\Sigma$ and $g(\Sigma)$ is the genus of $\Sigma$. It is easy to see that this means on any compact manifold there are only a finite number of cohomology classes that can be realized as the Euler class of a tight contact structure. This is in sharp contrast to the case of overtwisted contact structures that can realize any (even) cohomology class as their Euler class. Another contrast with overtwisted contact structures comes in work of Lisca and Mati\'c \cite{LiscaMatic97} who showed that there are homology spheres that admit an arbitrarily large (but finite) number of tight contact structures in a fixed homotopy class! A coarse understanding of tight contact manifolds was achieved in 2003 by Colin, Giroux, and Honda \cite{ColinGirouxHonda03,ColinGirouxHonda09} who showed that on any closed, oriented 3-manifold there are only finitely many homotopy classes of plane fields that contain a tight contact structure and if in addition the manifold is atoriodal, then the manifold has a finite number of tight contact structures! (We note that if the manifold has incompressible tori, one can always construct infinitely many tight contact structures using ``Giroux torsion".)

We finally turn to the classification of tight contact structures. Here our knowledge is significantly less than for overtwisted contact structures, but it does point to the classification being quite intricate and subtle. The first classification is, of course, due to Yasha. As noted above, in 1992 Yasha showed that $S^3$ admits a unique tight contact structure up to isotopy \cite{Eliashberg92a}. In 2000, the author showed that any lens space admitted a finite number of tight contact structures and proved that for one cohomology class on any lens space, there was a unique tight contact structure realizing it. But the next big breakthrough was due to Giroux \cite{Giroux00} and Honda \cite{Honda00a, Honda00b} who classified all tight contact structures on lens spaces, circle bundles over surfaces, and torus bundles. Their work comes down to studying ``convex surfaces", which is now a key tool in the study of contact structures. Convex surface theory grew out of Giroux's \cite{Giroux91} analysis of the notation of a convex contact manifold introduced by Yahsa and Gromov \cite{EliashbergGromov91}. 

Over the last 20 years, there has been a steady stream of classification results; for example, \cite{Ghiggini08, GhigginiLiscaStipsicz06, GhigginiLiscaStipsicz07, Tosun20, Wu06} classified tight contact structures on many (but not all!) small Seifert fibered spaces and \cite{Ghiggini05a} classified such structures on Seifert fibered spaces over $T^2$ with a single singular fiber. In addition, there have been several classifications of tight contact structures on some hyperbolic manifolds \cite{ConwayMin20, MinNonino23pre}. 

We did not state the specific classification results as they can be quite complicated, but the main point is that we only have classification results for quite special manifolds, and we see a rich and beautiful theory but still do not have a real understanding of the subtleties. Such as, why do some manifolds have large numbers of tight contact structures in a fixed homotopy type of plane field, while others do not? Why do some contact structures stay tight when pulled back to the universal cover and others do not? We clearly have a long way to go to obtain a good picture of the number and types of tight contact structures a given $3$-manifold supports.

As discussed above, one of the main ways to construct tight contact structures is to construct contact structures that have symplectic fillings, and at this point in the story, one might think that tightness of a contact structure might be equivalent to symplectically fillability. It turns out this is not the case as was first observed in 2002 by Honda and the author in \cite{EtnyreHonda02b}. Since then, there have been constructions of tight but not fillable contact structures on some circle bundles \cite{LiscaStipsicz03, LiscaStipsicz04} (but these manifolds admit other contact structures that are symplectically fillable) and on hyperbolic manifolds \cite{KalotiTosun17, LiLiu19} (these manifold admit no fillable contact structures at all), but overall, constructing tight but non-fillable contact structures seems difficult. 

%%%%%%%%%%%%%%%%%%%%%%%%%%%%%%%%%%%%
\subsection{Applications of tight contact structures}
%%%%%%%%%%%%%%%%%%%%%%%%%%%%%%%%%%%%
We have already discussed that the Bennequin inequality holds in any tight contact $3$-manifold. Thus we see that tight contact structures can ``see" subtle information about topology, that is they can give bounds on the minimal genus of a Seifert surface. But there are many other applications of tight contact structures!

We start with an impressive result of Yasha's. In \cite{Eliashberg92a}, Yasha used his classification of tight contact structures on $S^3$ together with a clever use of pseudo-holomorphic curves to reprove Cerf's theorem which states that any diffeomorphism of $S^3$ extends over the $4$-ball. Cerf's proof \cite{Cerf68} of his theorem was impressive and quite involved, while Yasha's proof was one page! (Of course, this is not really quite true, it takes quite a few pages to classify tight contact structures on $S^3$ and then one needs Gromov's theory of pseudo-holomorphic curves \cite{Gromov85}. So the proof is not really shorter, but what I find impressive is the two big pieces of Yasha's proof were developed for their own reasons and when combined, one gets a beautiful proof of this beautiful result! It really seems like contact geometry is trying to tell us deep things about topology.)

Another interesting application of contact geometry to topology is Ghiggini's proof that knot Heegaard Floer homology detects genus one fibered knots \cite{Ghiggini2008}. Specifically, given a non-fibered (genus one) knot $K$, Ghiggini constructs two tight contact structures on a manifold obtained by $0$ surgery on $K$, whose Heegaard Floer contact classes are different, but in the same grading. The idea behind this proof was generalized by Ni to prove that all fibered knots are detected by knot Heegaard Floer homology \cite{Ni07}. 

We next observe that Yasha and Thurston's result discussed above has major implications in the topology of $3$-manifolds. For example, among other results, it was essential in Kronheimer and Mrovka's proof of the Property P conjecture \cite{KronheimerMrowka04}, that says non-trivial surgery on a non-trivial knot in $S^3$ is not simply connected, and Ozsv\'ath and Szab\'o's proof that the Heegaard-Floer homology detects the Thurston norm of a manifold and the minimal Seifert genus of a knot \cite{OzsvathSzabo04a}.

We briefly sketch the proof of one of these results. Recall knowing the Thurston norm of a homology class is essentially the same as knowing the minimal genus of a surface representing that homology class (the norm is expressed in terms of the Euler characteristic of the surface, and care must be taken when the surface is a sphere). Suppose $\Sigma$ is a minimal genus surface in $M$ representing a homology class $h$. In \cite{Gabai83} Gabai tells us that there is a taut foliation $\mathcal{F}$ on $M$ that contains $\Sigma$ as a leaf. Yasha and Thurston's results above then say there is a symplectic structure on $X=M\times[-1,1]$ that symplectically fills $(M\times\{1\}, \xi_+)\cup (-M\times\{-1\},\xi_-)$ where the $\xi_\pm$ are positive and negative approximations of $\mathcal{F}$. Work of Yasha \cite{Eliashberg04} and the author \cite{Etnyre04a} then construct a closed symplectic manifold containing $X$. A non-vanishing result for the Heegaard-Floer invariant of a symplectic manifold in the spin$^c$ structure associated to the symplectic structure implies the non-vanishing of the Heegaard-Floer group $\HFplus(M,\spins_{\xi_+})$, where $\spins_{\xi_+}$ is the spin$^c$ structure associated to the contact structure $\xi_+$. We now recall that for any spin$^c$ structure $\spins$ with $\HFplus(M,\spins)\not= 0$, we have an adjunction inequality 
\[
|\langle c_1(\spins),[\Sigma]\rangle| \leq 2g(\Sigma)-2
\]
where $\langle c_1(\spins),[\Sigma]\rangle$ is the pairing of the Chern class of $\spins$ with the homology class of $\Sigma$ and $g(\Sigma)$ is the genus of $\Sigma$. Since $\Sigma$ was a leaf of $\mathcal{F}$ and $\xi_+$ is $C^0$ close to $\mathcal{F}$ it is easy to see that $\langle c_1(\spins_{\xi_+}),[\Sigma]\rangle=2g(\Sigma)+2$. Thus we see that we can detect the minimal genus of a non-trivial homology class in $M$ by seeing how $c_1$ of all the spin$^c$ structures on $M$ with non-zero Heegaard-Floer homologies evaluate on the homology class!

%%%%%%%%%%%%%%%%%%%%%%%%%%%%%%%%%%%%
\section{Higher dimensional overtwisted contact structures}
%%%%%%%%%%%%%%%%%%%%%%%%%%%%%%%%%%%%
As in dimension $3$, our discussion of the tight versus overtwisted dichotomy in higher dimensions begins with the questions of which manifolds admit contact structures. We saw in dimension $3$ that any oriented $3$-manifold admits such a structure, but in higher dimensions, there are obstructions to admitting a contact structure. Indeed, consider a co-oriented contact structure $\xi$ on a $(2n+1)$-manifold $M$. By choosing a contact form $\alpha$ for $\xi$ we obtain a symplectic structure $d\alpha$ on $\xi$. An almost contact structure on $M$ is a pair $(\alpha, \omega)$ where $\alpha$ is a non-degenerate $1$-form on $M$ and $\omega$ is a $2$-form that is symplectic when restricted to $\ker \alpha$. Clearly, if $M$ has a contact structure it has an almost contact structure. We also note that whether or not $M$ admits an almost contact structure is easy to determine; it is equivalent to the structure group of $TM$ reducing from $SO(2n+1)$ to $U(n)\oplus \mathbb{I}$ and this in turn can be determined by obstruction theory. So the real existence question in higher dimensions is whether or not any almost contact structure on a manifold $M$ is homotopic to a contact structure. 

There is a long history of studying this existence question and trying to define overtwisted contact structures in higher dimensions. In 1991, Geiges \cite{Geiges1991} showed that every almost contact structure on a simply connected $5$-manifold is homotopic to a contact structure, and in 1993, the case of ``highly connected" manifolds in all dimensions was addressed by Geiges \cite{Geiges1993}. In 1998, Geiges and Thomas \cite{GeigesThomas1998} gave a partial answer for $5$-manifolds with $\pi_1=\Z/2\Z$.  Despite this progress, the existence question remained elusive. For example, up until 2002, it was not known if tori $T^{2n+1}$ for $n>2$ admitted contact structures! (In 1979, Lutz \cite{Lutz1979} showed that $T^5$ admitted a contact structure, but for over 20 years nothing was known about $T^7$.) In 2002, Bourgeois \cite{Bourgeois02} showed that if $M$ admits a contact structure then so does $M\times T^2$. The first general result was in 2015 when Casals, Pancholi, and Presas \cite{CasalsPancholiPresas2015} showed that any almost contact structure on a $5$-manifold is homotopic to a contact structure.

Shortly after \cite{CasalsPancholiPresas2015}, Yasha, in collaboration with Borman and Murphy,  defined what it means for a contact structure to be overtwisted in all dimensions \cite{BormanEliashbergMurphy15}. The definition is a bit more complicated than in dimension $3$, so we refer to the paper for the details, but essentially a contact structure $\xi$ on $M^{2n+1}$ is overtwisted if there is an embedding of a piecewise smooth $2n$-disk that defines the germ of a model contact structure. With this definition in hand, they were able to show that any almost contact structure on a closed manifold is homotopic to an overtwisted contact structure that is unique up to isotopy. Thus the difficult existence problem is now solved in all dimensions!

Just as in dimension $3$, considerably more is true; they actually show that the space of contact structures containing a fixed overtwisted $D^{2n}$ and the space of almost contact structures also continuing this $D^{2n}$ are weakly homotopy equivalent. 

Prior to the above work on higher dimensional overtwisted contact structures, there was a previous attempt to generalize overtwistedness to higher dimensions. In \cite{Niederkrueger06} Neiderkr\'uger defined the notion of a ``plastikstufe". This was something like a parameterized version of an overtwisted disk, and he was able to show that having such an object implied that the contact structure could not be symplectically filled, just like for $3$-dimensional overtwisted contact structures! It is not hard to see that an overtwisted manifold contains plastikstufe, and in fact, if a contact manifold has a special type of plastikstufe it will be overtwisted \cite{CasalsMurphyPresas2019}.

Now that we have a notion of overtwisted in higher dimensions, we will define a tight contact structure to be one that is not overtwisted. As noted above if a contact manifold is symplectically fillable then it is not overtwisted, but one can wonder, as in dimension $3$, how tightness is related to symplectic fillability. The first progress on this was due to Massot,  Niederkr\"uger, and Wendl \cite{MassotNiederkrugerWendl13} who showed that there exist tight manifolds that are not symplectically fillable. There have been several similar results, but very recently Bowden, Gironella, Moreno, and Zhou, \cite{BowdenGironellaMorenoZhou24Pre} showed that tight but not symplectically fillable contact structures are everywhere. More explicitly, if $M$ is a contact manifold of dimension at least $7$ that admits a Stein fillable contact structure, then $M$ admits tight contact structures in the same homotopy class of almost contact structure that is not strongly fillable. The same is true for a $ 5$-manifold if the contact structure has first Chern class equal to zero. 

%and hence we may consider $\xi$ as a complex bundle by choosing a complex structure $J$ compatible with $d\alpha|_\xi$ (one may easily check that different choices for the complex structure compatible with $\alpha$ or other contact forms for $\xi$ are all homotopic as complex structures). Thus the tangent bundle to $M$ splits as $\xi\oplus \R$ where $\xi$ is a complex bundle. In other words, the structure group of $TM$ reduces from $SO(2n+1)$ to $U(n)\oplus \mathbb{1}$. 
%%%%%%%%%%%%%%%%%%%%%%%%%%%%%%%%%%%%
\section{Final thoughts}
%%%%%%%%%%%%%%%%%%%%%%%%%%%%%%%%%%%%
The tight versus overtwisted dichotomy was born just over 30 years ago in several seminal papers of Yasha. This dichotomy crystallized several insights of work from the past several years and provided a framework for the future of contact geometry. This framework has been the driving force for much of the work in the field and will certainly be so moving into the future.

% references
\bibliography{references}
\bibliographystyle{plain}
\end{document}